\documentclass{amsart}
\usepackage{amsthm,amssymb,booktabs}

\theoremstyle{plain}
\newtheorem{theorem}{Theorem}[section]
\newtheorem{lemma}[theorem]{Lemma}
\newtheorem{corollary}[theorem]{Corollary}
\newtheorem{prop}[theorem]{Proposition}
\theoremstyle{definition}

\newtheorem{example}[theorem]{Example}

\theoremstyle{remark}
\newtheorem{remark}[theorem]{Remark}

\newcommand{\NN}{{\mathbb N}}
\newcommand{\ZZ}{{\mathbb Z}}
\newcommand{\CC}{{\mathbb C}}
\newcommand{\fg}{{\mathfrak g}}
\newcommand{\hh}{{\mathfrak h}}
\newcommand{\fh}{{\mathfrak h}}
\newcommand{\Deltac}{\stackrel{\circ}{\Delta}}
\newcommand{\Lc}{\stackrel{\circ}{L}}
\newcommand{\hr}{\specialrule{\heavyrulewidth}{0pt}{0pt}}

\begin{document}

\title{Regular subalgebras of affine Kac--Moody algebras}
\author{Anna Felikson}
\address{Independent University of Moscow, B. Vlassievskii 11, 119002 Moscow, Russia}
\curraddr{Department of Mathematics, University of Fribourg, P\'erolles, Chemin du Mus\'ee 23, CH-1700 Fribourg, Switzerland} 
\email{felikson@mccme.ru}
\author{Alexander Retakh} 
\address{Department of Mathematics, Stony Brook University, Stony Brook, NY 11790, USA}
\email{retakh@math.sunysb.edu}
\author{Pavel Tumarkin}
\address{Independent University of Moscow, B. Vlassievskii 11, 119002 Moscow, Russia}
\curraddr{Department of Mathematics, Michigan State University, East Lansing, MI 48824, USA}
\email{pasha@mccme.ru}

\begin{abstract}
We classify regular subalgebras of affine Kac--Moody algebras in terms of their root systems. 
In the process, we establish that a root system of a subalgebra is always an intersection
of the root system of the algebra with a sublattice of its root lattice.

We also discuss applications to investigations of regular subalgebras of hyperbolic Kac--Moody algebras and conformally invariant subalgebras of affine Kac--Moody algebras.  In particular, we provide explicit formulae for determining all Virasoro charges in coset constructions that involve regular subalgebras.
\end{abstract}

\maketitle

\section{Introduction}
The main goal of this paper is to give a complete description of a large class of subalgebras of affine Kac--Moody algebras.  

Recall that both untwisted and twisted Kac--Moody algebras afford a uniform description in terms of a combinatorial datum: roughly speaking, an algebra $\fg$ is formed by a Cartan subalgebra $\hh$ and root subspaces indexed by a
discrete set $\Delta$, called a root system.  Our focus is on subalgebras whose root systems agree with that of $\fg$.

More precisely, consider a Kac--Moody algebra $\fg$.  A subalgebra $\fg_1\subset \fg$ is called {\it regular} if $\fg_1$ is invariant with respect to some Cartan subalgebra $\hh$ of $\fg$.  In other words, $\fg_1$ is a direct sum of a subspace of $\hh$ and subspaces of root spaces of $\fg$ (with respect to $\hh$).  In terms of root systems, $\fg_1$ is regular if its root system $\Delta_1$ is a specific subset of a root system $\Delta$ of $\fg$.

Our main result is the complete classification of regular affine subalgebras in both the untwisted and twisted case, see Theorems~\ref{th-ind}, \ref{th-dec}, and \ref{lower}.  We also describe a procedure for listing all regular subalgebras of affine Kac--Moody algebras.  As a by-product, we give a classification of those regular subalgebras of hyperbolic Kac--Moody algebras whose generalized Cartan matrices are positive definite or semidefinite.

Our starting point is the classification of affine Kac--Moody algebras presented in~\cite[Tables~Aff1--Aff3]{Kac} in terms of sets of simple roots of affine root systems.  For basic facts concerning Kac--Moody algebras, we also refer the reader to~\cite{Kac}.  In particular, we make use of the detailed description of affine root systems in~\cite[Proposition~6.3]{Kac}. 

Another important tool is the consideration of regular subalgebras from a geometric point of view. If $\fg_1\subset \fg$ is a regular subalgebra, then the natural embedding of respective root systems $\Delta_1\subset\Delta$ induces a natural embedding of corresponding Weyl groups $W_1\subset W$.  Weyl groups of affine Kac--Moody algebras are exactly affine reflection groups or, in other words, discrete groups of isometries of a Euclidean space generated by reflections~\cite[Proposition 3.13]{Kac}.  The relation with Weyl groups allows us, in turn, to use results of~\cite{eucl}, where all reflection subgroups of affine reflection groups were classified.        

Subalgebras of Kac--Moody algebras appear in various contexts.  In particular, just as an untwisted affine Kac--Moody algebra gives rise to a family of representations of the Virasoro algebra $Vir$ (via the Sugawara construction), every pair subalgebra/algebra also gives rise to a family of representations of $Vir$.  This is known as the coset construction.  Its applications include, for example, the explicit construction of all irreducible unitary highest weight representations of $Vir$ with charges between $0$ and $1$ \cite{GKO2}.  Another application is the construction of consistent conformally invariant theories: such subalgebras allow string compactification while preserving conformal covariance.  In this paper we restrict our attention to regular subalgebras only; however, in this case we are able to provide a recursive algorithm for determining all charges of $Vir$ obtained via the coset construction.  Since there is a formal analogue of the coset construction for twisted affine algebras (see e.g. \cite{KW}), we work with both untwisted and twisted algebras.

We note that our approach is not universal.  Not every subalgebra $\fg_1$ of an affine Kac--Moody algebra $\fg$ is amenable to our (primarily combinatorial) methods: to be regular $\fg_1$ must be semisimple (e.g. this rules out the Heisenberg subalgebra) and its root system must agree with that of $\fg$ (this rules out certain embeddings obtained via representations of $\fg_1$).  Nonetheless, we are able to deal effectively with a large and important class of subalgebras. 

The paper is organized as follows. In Section~\ref{roots} we list essential facts concerning root systems and introduce notations used throughout this paper.  In Sections~\ref{sub} and~\ref{subl} we prove a necessary and sufficient condition for a pair of root systems $\Delta_1\subset\Delta$ to define a pair algebra-subalgebra in terms of root lattices.  Section~\ref{af} contains the classification of  affine regular subalgebras of affine Kac--Moody algebras in terms of maximal subalgebras.  This includes both affine and non-affine subalgebras.  We simultaneously deal with both untwisted and twisted algebras.  In Section~\ref{coset}, we apply our results to the investigation of Virasoro charges obtained via the coset construction.  We also describe all maximal conformally invariant regular subalgebras.  Section~\ref{hyp} is devoted to subalgebras of hyperbolic Kac--Moody algebras.

\section{Regular subalgebras and root subsystems}
\label{roots}

Denote by $\Delta$ a root system of an affine Kac--Moody algebra $\fg$.  Our notations for the types of affine algebras follow the convention of~\cite[Tables Aff~1--Aff~3]{Kac}.  We use the same notation for the types of corresponding root systems.  Weyl groups are denoted as $\widetilde{X}_n$; we follow the convention of~\cite{V}.

For groups and algebras of small ranks we will also use the following notation:
\begin{align*}
 &\widetilde D_2=2\widetilde A_1 & &\widetilde B_1=\widetilde  A_1 &  &\widetilde C_1=\widetilde A_1 & 
&\widetilde D_3=\widetilde A_3 & &\widetilde B_2=\widetilde  C_2 & \\
 &D_2^{(1)}=2A_1^{(1)} & &B_1^{(1)}=\dot A_1^{(1)} & &C_1^{(1)}=\ddot A_1^{(1)} &
&D_3^{(1)}=A_3^{(1)} & &B_2^{(1)}=\dot C_2^{(1)} & &A_1^{(2)}=\ddot A_1^{(1)}
\end{align*}

Here one dot indicates that the shortest root of the corresponding root system has length $1$, while two dots indicate that the shortest root of the corresponding root system has length $\sqrt{2}$.  

Recall that root systems of {\it untwisted} algebras are listed in Table Aff~1; algebras of other types are called {\it twisted}. 

By abuse of terminology, we will call a root system of an affine (respectively, hyperbolic) Kac--Moody algebra an {\it affine (hyperbolic) root system}.   
A root system is called {\it decomposable} if it is a union of two mutually orthogonal root systems. Otherwise it is called {\it indecomposable}.  Indecomposable root systems correspond to simple algebras.

For a root system $\Delta$ we denote by $\Pi$ a set of simple roots of $\Delta$.
We keep standard notations $\Delta^{\mbox{re}}$ and $\Delta^{\mbox{im}}$ for real and imaginary roots of $\Delta$.
By $\Delta_+$ we mean the set of positive roots of $\Delta$ with respect to some fixed set of simple roots $\Pi$.

\medskip

Let $\Delta$ be a root system of a Kac--Moody algebra $\fg$.  Every regular subalgebra $\fg_1$ of $\fg$ has a root system $\Delta_1\subset\Delta$ that is closed with respect to addition, i.e. satisfies the following condition:
$$
\mbox{if\ }\alpha,\beta\in\Delta_1 \mbox{\ and\ } \alpha+\beta\in\Delta, \mbox{\ then\ } \alpha+\beta\in\Delta_1.
$$
A root system $\Delta_1\subset\Delta$ satisfying the condition above is called a {\it root subsystem} of $\Delta$. (Another term used in literature is subroot system.)  It  is easy to see that any root subsystem $\Delta_1\subset\Delta$ is a root system of a regular subalgebra of $\fg$.  Therefore, classifying regular subalgebras of $\fg$ is equivalent to classifying root subsystems of $\Delta$.
     
\begin{remark}  We wish to emphasize the difference between $\Delta_1$ and subsets 
of $\Delta$ used to construct representations of $\fg$.  The latter
are usually taken to contain a half of $\Delta$.  More precisely, the first
step in constructing representations of $\fg$ is to take a set $P\subset\Delta$
closed with respect to addition such that $P\cup -P=\Delta$ (see,
e.g. \cite{JK, F}).  Such a set $P$ corresponds to a Borel or, more generally,
parabolic subalgebra of $\fg$ whose representations yield those of $\fg$.
On the other hand, we are looking for semisimple subalgebras of $\fg$.  Thus
$\Delta_1$ is a root system itself.  For instance, unlike $P$, $\Delta_1=-\Delta_1$.
\end{remark}
     
To simplify checking if $\Delta_1$ is a root subsystem of a root system $\Delta$, we will restrict our attention to simple roots only.  Specifically, we have the following explicit criterion for $\Delta_1$ to be a root subsystem. Let $\Pi_1$ be a set of simple roots of $\Delta_1$. Then $\Delta_1$ is a root subsystem of $\Delta$ if and only if
\begin{equation}\label{eq.cond}
 \alpha-\beta\notin\Delta \mbox{\ for all\  } \alpha,\beta\in \Pi_1.
\end{equation}
This criterion is an immediate corollary of the following Proposition~\ref{fn}. 

\begin{prop}[{\cite[Theorem~3.1]{FN}}]
\label{fn}
Let $\Delta$ be a root system of a Kac--Moody algebra $\fg$.
Let $\beta_1,\dots,\beta_k\in \Delta^{\mbox{re}}$ be positive real roots such that
$\beta_i-\beta_j \notin \Delta$ for all $1\le i<j\le k$.
Let $\Delta'\subset \Delta$ be a minimal root system containing $\beta_1,\dots,\beta_k$.
Then $\Delta'$ is a root subsystem of $\Delta$
(and $\Delta'$ is a root system of a regular subalgebra of $\fg$). 
\end{prop}

We say that a root subsystem $\Delta_1\subset \Delta$ is {\it maximal} if for any root system $\Delta_2\subset \Delta$ containing $\Delta_1$ the system $\Delta_2$ is not a subsystem of $\Delta$. Clearly, maximal root subsystems correspond to maximal regular subalgebras.

\section{Root subsystems in root systems.}
\label{sub}

In this section we provide several technical statements that will drastically reduce the
computations necessary for classification.

We start with two elementary lemmas.

\begin{lemma}
\label{tower}
Let $\Delta_2\subset \Delta_1 \subset \Delta$ be root systems.

(1)
If $\Delta_1$ is a root subsystem of $\Delta$, and $\Delta_2$ is a root subsystem of $\Delta_1$, then $\Delta_2$ is a root subsystem of $\Delta$.

(2)
If $\Delta_2\subset \Delta_1$ is not a root subsystem,
then  $\Delta_2\subset \Delta$ is also not a root subsystem.

\end{lemma}

\begin{proof}

(1)
Let $\alpha,\beta\in \Delta_2$ and $\alpha+\beta\in \Delta$. To show that $\Delta_2$ is a root subsystem of $\Delta$, we check the definition by verifying that $\alpha+\beta\in \Delta_2$. Since $\Delta_2\subset \Delta_1$, it follows that $\alpha,\beta\in \Delta_1$.
By assumption, $\Delta_1\subset \Delta$ is a root subsystem, 
which implies that $\alpha+\beta\in \Delta_1$.
Since $\Delta_2\subset \Delta_1$, we also have  $\alpha+\beta\in \Delta_2$, so $\Delta_2\subset \Delta$ is a root subsystem.

(2)
By assumption there exist $\alpha,\beta\in \Delta_2$ such that $\alpha+\beta\in \Delta_1$ but $\alpha+\beta\notin \Delta_2$.
Since $\Delta_1\subset \Delta$, it follows that $\alpha+\beta\in \Delta$.  Hence $\Delta_2$ is not a root subsystem of $\Delta$.

\end{proof}

\begin{lemma}
\label{simplylaced}
Let $\Delta$ and $\Delta_1\subset \Delta$ be root systems, and let $\Delta$ be simply-laced (i.e. all roots in $\Delta$ are of the same length).
Then $\Delta_1$ is a root subsystem of $\Delta$.

\end{lemma}

\begin{proof}
Let $\alpha,\beta\in \Pi_1$ be such that $\alpha-\beta\in \Delta$.
Since $\alpha$ and $\beta$ are simple roots, 
the angle formed by  $\alpha$ and $\beta$
is not acute, so $\alpha-\beta$ is longer than $\alpha$ or $\beta$, which is impossible.
\end{proof}

\begin{remark} The same argument as in the proof of Lemma~\ref{simplylaced} also works in a non-simply-laced case.
To check condition~(\ref{eq.cond}) in a non-simply-laced case,
it is enough to check it for short roots $\alpha$ and $\beta$ (and the roots of middle length in the case
of $A_{2n}^{(2)}$).
\end{remark}

Lemma~\ref{k^n} describes a procedure for classifying root subsystems $\Delta_1\subset \Delta$ whose type coincides with the type of $\Delta$.  

\begin{lemma}
\label{k^n}
Let $\Delta$ and $\Delta_1\subset \Delta$ be root systems of affine algebras $\fg$ and $\fg_1$ of the same type and the same rank $n$. Let $W$ and $W_1\subset W$ be the corresponding Weyl groups.
%Let $n$ be a rank of $\Delta$ and $\Delta_1$. 
Then $[W:W_1]=k^n$ for some $k\in \ZZ_+$. Moreover,

(1) If $\fg$ is untwisted, then $\Delta_1\subset \Delta$ is a root subsystem for any $k$.   

(2) If $\fg$ is a twisted algebra of a type different from $D_{4}^{(3)}$,
 then $\Delta_1\subset \Delta$ is a root subsystem if and only if $k$ is odd.
 
(3) If $\Delta$ is of type $D_{4}^{(3)}$, then 
$\Delta_1\subset \Delta$ is a root subsystem if and only if $k$ is not a multiple of 3.
%(and not a subsystem otherwise).\\

\end{lemma} 
 
\begin{proof}
It follows from~\cite[Lemmas~9-11]{eucl} that if $W$ and $W_1\subset W$ are Euclidean reflection groups of the same type,
then either $[W:W_1]=k^n$ for some $k\in \ZZ_+$ or 
$W$ is of one of the types $\widetilde C_2$, $\widetilde G_2$, and $\widetilde F_4$.
In the latter case there exists a group $W_2$ of the same type as
$W$ such that
$W_1\subset W_2\subset W$, $[W:W_2]=k^n$ for some $k$, and $[W_2:W_1]=2,3$, or $4$ respectively
(depending on the type of $W$).
Checking the groups  $\widetilde C_2$, $\widetilde G_2$, and $\widetilde F_4$ directly, we see
that the root systems corresponding to the groups $W_2$ and $W_1$ are of different types
(their Dynkin diagrams differ by directions of arrows on the edges); however, the root systems corresponding to the groups $W_2$ and $W$ have the same type. Therefore,
the same type of $\Delta$ and $\Delta_1$ implies  $[W:W_1]=k^n$ for some $k\in \ZZ_+$. 

To determine which  embeddings of root systems are root subsystems,
we check condition~(\ref{eq.cond}).
A subgroup $W_1$ is determined (up to an automorphism of $W$) by the type of $W$ and the index   $[W:W_1]=k^n$,
so for each possible index we need to check one embedding only. 

We will use an explicit presentation of affine root systems given in~\cite[Prop.~6.3]{Kac}. Denote by $\Deltac$ the underlying finite root system of $\Delta$ (see~\cite[Chapter~6]{Kac}), and let $\Deltac_s$ and $\Deltac_l$ be the sets of short and long roots of $\Deltac$ respectively. 

(1)
By~\cite[Prop.~6.3]{Kac}, $\Delta=\{m\delta, \alpha+m\delta \, |\, \alpha\in \Deltac, m\in \ZZ \}$.
Then $\Delta_1=\{mk\delta, \alpha+mk\delta \, |\, \alpha\in \Deltac, m\in \ZZ \}$ where $k\in \NN$.
This implies that if $\alpha,\beta\in \Pi_1$, $\alpha-\beta\in \Delta$, then 
$\alpha-\beta\in \Delta_1$, which is impossible for simple roots $\alpha$ and $\beta$.

(2) 
If $\Delta\ne A_{2n}^{(2)}$, then~\cite[Prop.~6.3]{Kac} implies $\Delta=\{m\delta, \alpha+m\delta, \beta+2m\delta \, |\, \alpha\in \Deltac_{s}, 
\beta \in \Deltac_l, m\in \ZZ \}$.
Then $\Delta_1=\{mk\delta, \alpha+mk\delta, \beta+2mk\delta \,|\, \alpha\in \Deltac_{s}, 
\beta \in \Deltac_l, m\in \ZZ \}$ for some $k\in \NN$.
Notice, that $\Delta$ is a root system of a twisted algebra different from $A_{2n}^{(2)}$, thus
$\Delta_1$ contains two mutually orthogonal short simple roots $\alpha_1$ and $-\beta_1+k\delta$
(where $\alpha_1,\beta_1\in \Deltac_1$ and $\alpha_1+\beta_1\in \Deltac_1$). 
So, if $k$ is even, we have $\alpha_1-(-\beta_1+k\delta)\in \Delta$, and $\Delta_1$ is not a root subsystem of $\Delta$.
If $k$ is odd, then we check all pairs of short simple roots $\alpha,\beta\in \Pi_1$
and never obtain $\alpha-\beta\in \Delta$. 

If $\Delta=A_{2n}^{(2)}$, then $\Pi$ contains only one short simple root. So, if  
$\alpha,\beta\in \Pi_1$ and $\alpha-\beta\in \Delta_1$, then $\alpha$ and $\beta$ are two orthogonal roots of
the middle length. In this case $\alpha-\beta$  never belongs to $\Delta$.
Therefore, each root system $\Delta_1\subset \Delta$ of type $A_{2n}^{(2)}$ is a root 
subsystem of $\Delta$. On the other hand, if $W_1\subset W$ are Weyl groups 
of type $\widetilde C_n$ such that $[W:W_1]=k^n$ for an even $k$,
then the root system $\Delta_1$ consists of roots of only two different lengths.
In this case $\Delta_1$ is not a root system of  type $A_{2n}^{(2)}$.
In case of odd $k$ $\Delta_1$ contains roots of all three lengths,
so  it is a root system of  type $A_{2n}^{(2)}$.

(3) The proof is similar to the proof of (2) for $\Delta\ne A_{2n}^{(2)}$.

\end{proof}

\section{Root subsystems and root lattices}
\label{subl}

The goal of this section is to prove Theorem~\ref{th-sublattice}. 

The following lemma holds for finite root systems as well as for affine and hyperbolic ones.
The finite case is well known and the hyperbolic case may be found in~\cite{roots}.
Here we prove the affine case. For completeness we include also a proof of the finite case;
see the claim below.

\begin{lemma}
\label{sublattice}
Let  $\Delta$ and $\Delta_1\subset \Delta$ be affine root systems
of the same rank.
Let $L$ and $L_1\subset L$ be the corresponding root lattices.
Then $\Delta_1$ is a root subsystem of $\Delta$ if and only if \quad 
$\Delta_1 = \Delta \cap L_1.$

\end{lemma}

\begin{proof}
First, assume that $\Delta_1 = \Delta \cap L_1$ but $\Delta_1$ is not a root subsystem
of $\Delta$.
Then there exist $\alpha,\beta \in \Pi_1$ such that
$\alpha-\beta \in \Delta$.  The condition $\Delta_1 = \Delta \cap L_1$ implies that 
$\alpha-\beta\in \Delta_1$.  This is impossible since 
$\alpha$ and $\beta$ are simple roots of $\Delta_1$.

Now suppose $\Delta_1$ be a root subsystem of $\Delta$.
Suppose also that $\Delta_1\ne \Delta\cap L_1$,
so that there exists $\alpha\subset \Delta$ such that $\alpha\in L_1\setminus \Delta_1$.

Denote by $\Deltac$ and $\Deltac_1$ the underlying finite root 
systems of $\Delta$ and $\Delta_1$ respectively (see~\cite[6.3]{Kac}), 
and suppose that $\Delta_1\ne A_{2n}^{(2)}$. Then $L$ is generated over $\ZZ$ 
by $\Deltac$ and 
$\delta$, and $L_1$ is generated by $\Deltac_1$ 
and $k\delta$ for some integer $k\ge 1$. Let $\alpha=\alpha_0+l\delta$, 
$\alpha_0\in\Deltac$. 
Clearly, $\Deltac_1$ is a root subsystem of $\Deltac$. 
Denote by $\Lc$ and $\Lc_1\subset \Lc$ 
the root lattices of $\Deltac$ and $\Deltac_1$,
respectively.

\noindent
{\bf Claim.} $\Deltac_1 = \Deltac \cap \Lc_1.$\\
We will prove this claim for the case when $\Deltac_1$ is a maximal root subsystem of $\Deltac$. This will immediately imply the claim in the general case.  

By description of maximal regular subalgebras of semisimple Lie algebras (see~\cite[Theorem~5.5]{Dyn}), in this case $\Lc_1$ is a proper sublattice of $\Lc$. The intersection $\widetilde\Delta_1 = \Deltac \cap \Lc_1$ is a root subsystem of $\Deltac$ containing $\Deltac_1$. Since $\Deltac_1\subset\Deltac$ is maximal, we see that $\widetilde\Delta_1 = \Deltac_1$.

\medskip 

Applying the claim above we have $\alpha_0\in\Deltac_1$. 
Moreover, $\alpha\in L_1$ implies $k\mid l$. 
Now, if $\Delta_1$ is a root system of an untwisted algebra, we see that 
$\alpha\in\Delta_1$, contradicting our assumption. 
If $\Delta_1$ is a root system of a twisted algebra and 
$\alpha\notin\Delta_1$, we see that $\alpha_0$ is a long root 
of $\Deltac_1$, and either 
$3\not\:\mid (l/k)$ (in the case of $\Delta_1=D_4^{(3)}$), 
or $l/k$ is odd (otherwise). In both 
cases, at least one of $\alpha+k\delta$ and 
$\alpha-k\delta$ is a root of $\Delta_1$. Combined with the assumption that 
$\alpha\in\Delta_1$ and $\pm k\delta\in\Delta_1$, this implies 
that $\Delta_1$ is not a root subsystem of $\Delta$. 
The contradiction proves the lemma for all $\Delta_1\ne A_{2n}^{(2)}$.         

In the case $\Delta_1= A_{2n}^{(2)}$, we have $\Delta= A_{2n}^{(2)}$. 
Then $\Deltac_1=\Deltac$.  Furthermore,
$L$ is generated over $\ZZ$ by $\Deltac$ and 
$(\alpha_n+\delta)/2$ (where $\alpha_n$ is a long root of $\Deltac$), 
and $L_1$ is generated by $\Deltac$ and $(\alpha_n+(2k-1)\delta)/2$ 
for some integer $k\ge 1$. 
Two cases are possible: either $\alpha=(\alpha_n+(2l-1)\delta)/2$ for a long root 
$\alpha_n\in\Deltac$ and an integer $l$, 
or $\alpha=\alpha_0+l\delta$ 
for $\alpha_0\in\Deltac$ and an integer $l$.

Let $\alpha=(\alpha_n+(2l-1)\delta)/2$. Since $\alpha\in L_1$, we have $(2k-1)\mid(2l-1)$. 
This implies that $\alpha\in\Delta_1$.

Now let $\alpha=\alpha_0+l\delta$. Again, $\alpha\in L_1$ implies $(2k-1)\mid(2l-1)$. 
Since $\alpha\notin\Delta_1$, we see that $\alpha_0$ is a long root of 
$\Deltac_1$, and $(2l-1)/(2k-1)$ is odd. 
Hence, $\alpha+(2k-1)\delta\in\Delta_1$. 
Furthermore, $(2k-1)\delta$ is also a root of $\Delta_1$, which contradicts the assumption that 
$\Delta_1$ is a root subsystem of $\Delta$. This completes the proof of the lemma. 
\end{proof}

\begin{remark}
\label{any}
Note that to prove sufficiency of condition $\Delta_1 = \Delta \cap L_1$, we do not need
to assume that $\Delta_1$ and $\Delta$ have the same rank.
\end{remark}

\begin{corollary}
\label{proper sublattice}
Let $\Delta$ and $\Delta_1\subset \Delta$ be affine root systems
of the same rank.
Let $L$ and $L_1\subset L$ be the corresponding root lattices,
and let $W$ and $W_1$ be the Weyl groups of $\Delta$ and $\Delta_1$.
Assume that $W_1$ is a maximal reflection subgroup of $W$.
Then $\Delta_1$ is a root subsystem of $\Delta$ if and only if 
$L_1$ is a proper sublattice of $L$.

Moreover, this sublattice has index $2$ unless
$(\Delta_1,\Delta)= (A_2^{(1)},G_2^{(1)})$ or $(G_2^{(1)},D_4^{(3)})$.
In the latter cases the index equals 3.

\end{corollary}

\begin{proof}

First, suppose that $L_1=L$. Then 
$\Delta\cap L_1=\Delta\cap L=\Delta\ne \Delta_1$. Lemma~\ref{sublattice} implies that 
 $\Delta_1$ is not a root subsystem of $\Delta$. 

Now, let $L_1$ be a proper sublattice of $L$.
By Lemma~\ref{sublattice} it suffices to check that $\Delta_1 = \Delta \cap L_1$.
For $\Delta^{im}_1$ the assertion is obvious.
Assume that there exists $\alpha \in \Delta^{re}$ such that $\alpha\in L_1$ 
such that $\alpha\notin \Delta_1$. 

Consider the
subgroup $G_0$ of $W$ generated by reflections with respect to all roots in $L_1$.
Obviously $W_1\subset G_0$. Since $W_1$ is maximal in $W$ 
and $G_0 \ne W_1$, it follows that
$G = W_1$. Hence each simple root $\alpha$ in $\Delta$ can be represented as a linear 
combination of roots in $L_1$ with integer coefficients.  Therefore $\alpha$ belongs 
to $L_1$ itself. Thus,
$L = L_1$ contradicting the assumption.

A direct calculation shows that for each but two cases, the
roots of the subsystem generate an index-two sublattice of the root lattice,
while in the remaining cases the index of the sublattice equals 3.
\end{proof}

To prove the main result of this section, we drop the assumption that $\Delta$ and its subsystem $\Delta_1$ have the same rank. 

\begin{theorem}
\label{th-sublattice}
Let  $\Delta$ and $\Delta_1\subset \Delta$ be affine root systems.
Let $L$ and $L_1\subset L$ be the corresponding root lattices.
Then $\Delta_1$ is a root subsystem of $\Delta$ if and only if \quad 
$\Delta_1 = \Delta \cap L_1.$

\end{theorem}

\begin{proof}
We need to prove one implication only (see Remark~\ref{any}). Namely, suppose that $\Delta_1$ is a root subsystem of $\Delta$. Then we need to show that $\Delta_1 = \Delta \cap L_1$. Moreover, we may assume that $\Delta_1\subset \Delta$ is a maximal root subsystem. 

If ranks of $\Delta$ and $\Delta_1$ are equal, then the statement follows from Lemma~\ref{sublattice}. Hence, we can assume that the rank of $\Delta_1$ is less than the rank of $\Delta$. In particular, $L_1$ is a proper sublattice of $L$. 

Consider the set $\Delta'=\Delta \cap L_1$. Clearly $\Delta'$ is a root subsystem of $\Delta$ containing $\Delta_1$. Since $L_1$ does not coincide with $L$, $\Delta'$ is a proper subsystem of $\Delta$. Now maximality of $\Delta_1\subset \Delta$ implies $\Delta'=\Delta_1$ and the theorem is proved.
\end{proof}

\section{Classification of regular subalgebras}
\label{af}

Given an affine root system $\Delta_1$ contained in an affine root system $\Delta$,  
we can use condition~(\ref{eq.cond}) to determine if $\Delta_1$ is a root substystem of $\Delta$.
Our goal is to list all affine root subsystems of root systems of every affine Kac--Moody algebra.

Checking condition~(\ref{eq.cond}) is straightforward, and Lemmas~\ref{tower},~\ref{simplylaced}, and~\ref{k^n} help to reduce the computations.
Thus essentially we only need to list all possible pairs $(\Delta, \Delta_1)$. This is done as follows.

For each of the  finitely many types of affine algebras (listed in~\cite[Tables~Aff~1-Aff~3]{Kac})
we consider its root system $\Delta$ and Weyl group $W$. We list all possibilities for its reflection subgroup 
$W_1\subset W$
using the results of~\cite{eucl}.
Then for each subgroup $W_1\subset W$ we check all possible root systems $\Delta_1\subset \Delta$ 
corresponding to such a subgroup.   For example, a subgroup of type $\widetilde B_n$ may correspond to
a root system of either type $B_n^{(1)}$ or type $A_{2n-1}^{(2)}$, while a subgroup of type
$\widetilde C_n$ may correspond to one of the types $C_n^{(1)}$, $D_{n+1}^{(2)}$, and $A_{2n}^{(2)}$.  
Furthermore, we consider not only the types of subgroups (and corresponding root systems)
but also all different embeddings of subgroups of this type.
Most embeddings of Euclidean reflection subgroups of finite index of type $\widetilde X$ in Euclidean reflection
groups of type $\widetilde Y$ can be distinguished by the indices of subgroups.
So, to indicate which embeddings correspond to subalgebras we specify the indices of respective subgroups
(in case of finite index).
If a subgroup is of infinite index, we check only if a root system $\Delta_1$ of a given type can be a root subsystem of $\Delta$.  (If the answer is affirmative, we place this subsystem in the resulting tables without further distinguishing it from other embeddings, which may not be root subsystems).
For some types of root systems  $\Delta$ and $\Delta_1\subset \Delta$, the system $\Delta_1$ may be either short or long
(for example, for $D_4^{(1)}\subset F_4^{(1)}$). In such cases we check all possibilities and
state in the tables below whether the root subsystem is short or long.  

\begin{remark}
It was stated in~\cite[Theorem~1]{eucl} that if $W$ is an indecomposable group 
of a given type $Y$, then all reflection subgroups of a given type
$X$ can be distinguished by their indices (modulo the automorphism group of $W$).
However, while preparing this paper, we discovered an exception to this rule. 

Namely, if $W$ and $W_1$ are of types $\widetilde F_4$ and $\widetilde C_4$ respectively,
 there are two different embeddings of a subgroup of type $\widetilde C_4$ into a group  
of type $\widetilde F_4$ with the same index $[W:W_1]= 24k^4$ for any $k\in \NN$.
Specifically, there is a chain of embeddings $\widetilde C_4\subset \widetilde B_4 \subset \widetilde F_4$,
where an automorphism of $\widetilde C_4$ may be extended to an automorphism of $\widetilde B_4$,
but not to an automorphism of $\widetilde F_4$.

The oversight resulted from an omission in a case by case check at the very end of the proof of Theorem~1. 
\end{remark}

To make our tables more readable, we list maximal root subsystems
only (others can be easily retrieved from the lists of maximal ones).
Note that a decomposable maximal root subsystem does not necessarily correspond to a maximal subgroup,
namely, for a chain of root systems $\Delta_1\subset \Delta_2 \subset \Delta_3$
it is possible that  $\Delta_2\subset \Delta_3$ is not a root subsystem, while 
$\Delta_1\subset \Delta_3$ is.

First, we consider indecomposable affine root subsystems of maximal rank,
so that $W_1\subset W$ is an indecomposable Euclidean reflection subgroup of finite index.
Then we allow decomposable $W_1$.  Further on, we move to affine root subsystems of lower rank
(and then $[W:W_1]$ is infinite).  Finally we describe remaining root subsystems
that correspond to direct products of a spherical and a Euclidean reflection groups.

\subsection{Indecomposable affine root subsystems of maximal rank}

Below we provide a typical procedure of classifying maximal rank root subsystems of the root system of an untwisted algebra.
Dealing with twisted algebras usually involves more computations but is not otherwise different. 

\begin{example}[{\bf Indecomposable affine root subsystems of $B_n^{(1)}$, $n\ge 3$}]
\label{mbu}
The Weyl group $W$ is of the type $\widetilde B_n$.  Results in~\cite{eucl} imply that 
its indecomposable finite index subgroups have types $\widetilde B_n$, $\widetilde C_n$,
and $\widetilde D_n$.  For subgroups of these types we should consider root systems of types
$B_n^{(1)}$ and $A_{2n-1}^{(2)}$ (for $\widetilde B_n$), $C_n^{(1)}$, $D_{n+1}^{(1)}$, and $A_{2n}^{(2)}$ (for $\widetilde C_n$), and $D_n^{(1)}$ (for $\widetilde D_n$).  The root system of type $A_{2n}^{(2)}$ is irrelevant because it contains roots of three different 
lengths, while $B_n^{(1)}$ only has roots of two different lengths. 
Furthermore, short roots of $B_n^{(1)}$ are mutually orthogonal,
which implies that neither  $A_{2n-1}^{(2)}$ nor $C_n^{(1)}$ can be contained in $B_n^{(1)}$.
So, we need to check only root systems of types $B_n^{(1)}$, $D_{n+1}^{(2)}$, and $D_n^{(1)}$.
Notice that  embeddings of groups $\widetilde B_n$, $\widetilde C_n$, and $\widetilde D_n$ into
$\widetilde B_n$ are determined by indices up to  isomorphisms of $\widetilde B_n$, so
it suffices to check one embedding for each possible index.
Explicit check shows that embeddings of $B_n^{(1)}$ and $D_n^{(1)}$ satisfy condition~(\ref{eq.cond}),
while that of $D_{n+1}^{(2)}$ does not. 
\end{example}

\begin{theorem}
\label{th-ind}
Let $\Delta$ be an indecomposable affine root system
and $\Delta_1\subset \Delta$ a maximal root subsystem.
If $\Delta_1$ is an indecomposable affine root subsystem of maximal rank,
then $\Delta_1\subset \Delta$ is one of the root subsystems listed in Table~\ref{indecomp}. 

\end{theorem}

\renewcommand\arraystretch{1.3}
\begin{table}[h]
\caption{\label{indecomp}Indecomposable maximal affine root subsystem of maximal rank.
The right column lists possible indices of Weyl subgroups $W_1\subset W$ corresponding to maximal 
root subsystems $\Delta_1\subset \Delta$
($p$ is any prime number).}

\begin{tabular}{ccc}

\begin{tabular}{l|l|l}
\hr
$\Delta$ & $\Delta_1$ & $[W:W_1]$  \\
\hr
$A_1^{(1)}$ & $A_1^{(1)}$ & $p$  \\
\hline
$A_n^{(1)}$ \ $(n\ge 2)$  & $A_n^{(1)}$ & $p^n$ \\
\hline
$B_n^{(1)}$ \ $(n\ge 3)$  & $B_n^{(1)}$ & $p^n$  \\
                          & $D_n^{(1)}$ & $2$  \\
\hline
$C_n^{(1)}$ \ $(n\ge 2)$  & $C_n^{(1)}$ & $p^n$  \\
\hline
$D_n^{(1)}$ \ $(n\ge 4)$  & $D_n^{(1)}$ & $p^n$  \\
\hline
$G_2^{(1)}$               & $\ddot A_2^{(1)}$ & $2$  \\
                          & $G_2^{(1)}$ & $p^2$  \\
\hline
$F_4^{(1)}$               & $B_4^{(1)}$ & $3$  \\
                          & $F_4^{(1)}$ & $p^4$  \\
\hline
$E_6^{(1)}$ & $E_6^{(1)}$ & $p^6$ \\
\hline
$E_7^{(1)}$ & $A_7^{(1)}$ & $2^4\cdot 3^2$  \\
            & $E_7^{(1)}$ & $p^7$ \\
\hline
$E_8^{(1)}$ & $A_8^{(1)}$ & $2^7 \cdot 3^2 \cdot 5$  \\
            & $D_8^{(1)}$ & $2\cdot 3^3\cdot 5$  \\
            & $E_8^{(1)}$ & $p^8$  \\

\multicolumn{3}{c}{}\\
\multicolumn{3}{c}{}

\end{tabular}

&\quad &
\begin{tabular}{l|l|l}
\hr
$\Delta$ & $\Delta_1$ & $[W:W_1]$  \\
\hr
$A_2^{(2)}$ &  $A_2^{(2)}$ & $p$, $p\ne 2$\\
            &  $\dot A_1^{(1)}$ & $2$ \\
            &  $\ddot A_1^{(1)}$ & $2$ \\
\hline
$A_{2n}^{(2)}$, $n\ge 2$  &  $A_{2n}^{(2)}$ & $p^n$, $p\ne 2$\\
  &  $A_{2n-1}^{(2)}$ & $2$\\
 &  $B_{n}^{(1)}$ & $2^{n+1}$\\
\hline

$A_4^{(2)}$ &  $D_3^{(2)}$ & $2$\\
 
\hline

$A_{2n-1}^{(2)}$ $n\ge 3$&  $A_{2n-1}^{(2)}$ & $p^n$, $p\ne 2$\\
  &  $C_{n}^{(1)}$ & $2^{n-1}$\\
\hline
$D_{n+1}^{(2)}$ $n\ge 2$&  $D_{n+1}^{(2)}$ & $p^n$, $p\ne 2$\\
  &  $B_{n}^{(1)}$ & $2$\\
\hline
$D_{3}^{(2)}$ &  $C_{2}^{(1)}$ & $2$\\
\hline

$E_6^{(2)}$ &  $E_6^{(2)}$ & $p^4$, $p\ne 2$\\
  &  $F_{4}^{(1)}$ & $4$\\
  &  $C_{4}^{(1)}$ & $2^3\cdot3$\\
\hline
$D_4^{(3)}$ &  $D_4^{(3)}$ & $p^2$, $p\ne 3$\\
  &  $G_{2}^{(1)}$ & $3$\\
  &  $\dot A_{2}^{(1)}$ & $2$\\

\end{tabular}

\end{tabular}
\medskip
\end{table}
\renewcommand\arraystretch{1}

\begin{remark}
Our notation in Tables~\ref{indecomp} and~\ref{decomp} is as follows. If a root system of some type $X_n^{(r)}$ may be embedded in $\Delta$ with different root lengths, we write $\dot X_n^{(r)}$ for a ``short'' embedding and $\ddot X_n^{(r)}$ for the ``long'' one.
\end{remark}

\subsection{Decomposable affine root subsystems of maximal rank}

\begin{example} [{\bf Decomposable affine root subsystems of $B_n^{(1)}$, $n\ge 3$}]
\label{mbd}
As shown in Example~\ref{mbu}, $B_n^{(1)}$ may contain root systems of types $B_n^{(1)}$, $D_n^{(1)}$,
and $D_{n+1}^{(2)}$
(not necessary as root subsystems).
Furthermore, it follows from~\cite[Table~5]{eucl} that any maximal subgroup of $\widetilde B_n$ 
is $\widetilde B_m +\widetilde B_{n-m}$ ($m\ge 1$), and similarly
maximal subgroups of $\widetilde C_n$ and  $\widetilde D_n$ are $\widetilde C_m +\widetilde C_{n-m}$
and $\widetilde D_m+\widetilde D_{n-m}$, respectively.
So, looking through the chains of maximal subgroups
(and taking into account that short roots of $B_n^{(1)}$ are mutually orthogonal) 
 we conclude that any root subsystem $\Delta_1$ of  $B_n^{(1)}$ splits into components of types $B_m^{(1)}$, $D_m^{(1)}$
and $D_{m+1}^{(2)}$ (with possibly different ranks of components). 
If $\Delta_1$ has at least one component of type $D_{m+1}^{(2)}$,
then $\Delta_1$ is not a root subsystem.  Indeed, if $\alpha$ and $\beta$ are short simple roots of $D_{m+1}^{(2)}$,
then $\alpha-\beta$ belongs to $B_n^{(1)}$. 
Similarly, if $\Delta_1$ is a root subsystem, then it contains at most one component of type $B_m^{(1)}$.
Therefore, any decomposable affine root subsystem (of maximal rank) in  $B_n^{(1)}$
can be written  as 
$\varepsilon B_{m_0}^{(1)}+D_{m_1}^{(1)}+\dots+D_{m_s}^{(1)}$, where $\varepsilon=0$ or $1$
and $\varepsilon m_0+m_1+\dots+m_s=n$ for $m_0\ge 1$ and $m_i\ge 2$, $i=1,\dots,s$.
Clearly, all these root systems are root subsystems of  $\widetilde B_n$, and the only maximal subsystems are $B_{m}^{(1)}+D_{n-m}^{(1)}$ for $m=1,\dots,n-2$.
\end{example}

Similar straightforward calculations for every indecomposable affine root system
imply the following theorem.

\begin{theorem}\label{th-dec}
Let $\Delta$ be an indecomposable affine root system
and $\Delta_1\subset \Delta$ its maximal root subsystem.
If $\Delta_1$ is a decomposable root subsystem of maximal rank,
then $\Delta_1\subset \Delta$ is one of the subsystems listed in Table~\ref{decomp}. 

\end{theorem}

\renewcommand\arraystretch{1.3}
\begin{table}[h]
\caption{\label{decomp}Decomposable maximal affine root subsystems of maximal rank.}
\begin{tabular}{l|l|l}
\hr
$\Delta$ & $\Delta_1$ & $[W:W_1]$  \\
\hr
$B_n^{(1)}$ & $B_m^{(1)}+D_{n-m}^{(1)}$ & $4{n\choose m}$  \\
\hline
$C_n^{(1)}$ & $C_m^{(1)}+C_{n-m}^{(1)}$ & ${n\choose m}$  \\
\hline
$D_n^{(1)}$ & $D_m^{(1)}+D_{n-m}^{(1)}$ & $4{n\choose m}$  \\
\hline

$E_6^{(1)}$& $A_5^{(1)}+A_1^{(1)}$& $2^3\cdot 3^2$\\
            & $3A_2^{(1)}$ & $2^4\cdot 3^2\cdot 5$\\ 
\hline

$E_7^{(1)}$& $D_6^{(1)}+A_1^{(1)}$& $2\cdot 3^2\cdot 7$\\
            & $A_5^{(1)}+A_2^{(1)}$& $2^5\cdot 3^2\cdot 7$\\ 
\hline

$E_8^{(1)}$& $E_7^{(1)}+A_1^{(1)}$& $2^4\cdot 3\cdot 5$\\
            & $E_6^{(1)}+A_2^{(1)}$& $2^6\cdot 3\cdot 5\cdot7$\\ 
            & $2A_4^{(1)}$     & $2^8\cdot 3^3\cdot 5\cdot7$\\   
\hline

$F_4^{(1)}$& $\dot A_2^{(1)}+\ddot  A_2^{(1)}$ & $2^5\cdot 3$ \\
            & $C_3^{(1)}+\ddot A_1^{(1)}$ & $2^7\cdot 3$\\
\hline
 $G_2^{(1)}$& $\dot A_1^{(1)}+\ddot  A_1^{(1)}$ & $6$ \\

 \hline           

$A_{2n-1}^{(2)}$ & $A_{2m-1}^{(2)}+A_{2n-2m-1}^{(2)}$ & $2{n\choose m}$  \\
\hline
$A_{2n}^{(2)}$ & $A_{2m}^{(2)}+A_{2n-2m-1}^{(2)}$ & $2{n\choose m}$  \\
               & $D_{m}^{(1)}+A_{2n-2m}^{(2)}$ & $2^{m+2}{n\choose m}$  \\
               & $\ddot D_{3}^{(2)}+A_{2n-4}^{(2)}$ & $2{n\choose m}$  \\

\hline
$D_{n+1}^{(2)}$ & $D_{n-m}^{(1)}+D_{m+1}^{(2)}$ & $4{n\choose m}$  \\

\hline
$E_6^{(2)}$& $A_5^{(2)}+\ddot A_1^{(1)}$& $2^3\cdot 3$\\
\hline

 $D_4^{(3)}$& $\dot A_1^{(1)}+\ddot  A_1^{(1)}$ & $6$ \\

\end{tabular}
\medskip
\end{table}
\renewcommand\arraystretch{1}

\subsection{Affine root subsystems of lower rank}

To classify root subsystems of lower rank we use the same procedure as in the case of maximal rank.
More precisely, we consider chains of subgroups $G_1\subset G_2\subset \dots\subset G_s=W$,
where $W$ is a Weyl group of $\fg$ and each $G_i$ a maximal subgroup of $G_{i+1}$.
Then we check which root systems correspond to $G_i$ and which of those are root subsystems.
Since we obtain root subsystems of ranks smaller than that of $\fg$,
one of the subgroups $G_i\subset G_{i+1}$ must be a maximal subgroup of infinite index.
Such subgroups are described in~\cite[Theorem~3]{eucl}.
Namely, these maximal subgroups are of the following types:
$\widetilde A_k+ \widetilde A_{n-1-k}\subset \widetilde  A_n$ ($k \le n-1$);  
$\,\widetilde D_{n-1}\subset \widetilde D_n$; 
$\,\widetilde A_{n-1} \subset \widetilde D_n$;  $\,\widetilde  D_5\subset \widetilde E_6$; and $\,\widetilde  E_6 \subset \widetilde  E_7$.

\begin{example} [{\bf Lower rank affine root subsystems of $B_n^{(1)}$, $n\ge 3$}]
As in Example~\ref{mbd} we see that 
any affine root subsystem of $B_n^{(1)}$ can be written as
$\varepsilon B_{m_0}^{(1)}+D_{m_1}^{(1)}+\dots+D_{m_s}^{(1)}$, where $\varepsilon=0$ or $1$.
If $\varepsilon m_0+m_1+\dots+m_s< n$, then this is a root subsystem of $B_{n-1}^{(1)}$,
which in its turn is a root subsystem of $B_{n}^{(1)}$.
So, the only maximal root subsystem of lower rank is $B_{n-1}^{(1)}$.

Subgroups described in~\cite[Theorem~3]{eucl} do not affect the list of maximal root subsystems
of $B_{n}^{(1)}$.
Indeed, among the groups $\widetilde A_n$, $\widetilde D_n$, $\widetilde E_6$, and $\widetilde E_7$, 
only $\widetilde D_n$ is a  subgroup of $\widetilde B_n$.
The corresponding root system $D_{n}^{(1)}$ is a root subsystem of $B_{n}^{(1)}$,
thus none of root subsystems of $D_{n}^{(1)}$ can  be maximal in $B_{n}^{(1)}$.

\end{example}

Applying the above procedure to each indecomposable affine root systems,
we obtain the following theorem.

\begin{theorem}
\label{lower}
Let $\Delta$ be an indecomposable affine root system
and $\Delta_1\subset \Delta$ a maximal affine root subsystem.
If the rank of $\Delta_1$ is less than the rank of $\Delta$, 
then $(\Delta,\Delta_1)$ is one of the following pairs:
\begin{center}
\begin{tabular}{cc}
\begin{tabular}{l}
$(A_n^{(1)},A_k^{(1)}+A_{n-1-k}^{(1)})$, \\
$(B_n^{(1)},B_{n-1}^{(1)})$, \\
$(C_n^{(1)},A_{n-1}^{(1)})$, \\
$(D_n^{(1)},A_{n-1}^{(1)})$,\\
\end{tabular}
&
\begin{tabular}{l} 
$(E_7^{(1)},E_{6}^{(1)})$, \\
$(E_6^{(1)},D_{5}^{(1)})$, \\
$(A_{2n-1}^{(2)},A_{n-1}^{(1)})$,\\ 
$(D_{n+1}^{(2)},D_{n}^{(2)})$. \\
\end{tabular}
\end{tabular}
\end{center}
\end{theorem}

\subsection{Non-affine root subsystems}

Let $\Delta$ be an affine root system. Consider a regular subsystem $\Delta_1\subset \Delta$. 
Clearly $\Delta_1$ is a direct sum of an affine and a finite root systems $\Delta_a$ and $\Delta_s$.
 Let
$\Delta_s=\Delta_1\cup\dots\cup\Delta_k$, where $\Delta_i$ are indecomposable and mutually orthogonal.
Let $\Pi_i$ be a set of simple roots of $\Delta_i$, $1\leq i\leq k$.
Denote by $\Pi_s=\Pi_1\cup \dots\cup\Pi_k$ a set of simple roots of $\Delta_s$.
Add to each $\Pi_i$ a root $\beta_i=\theta_i+k_i\delta$, where $\theta_i$ is the lowest root of $\Delta_i$
and $k_i$ is the least positive integer such that $\theta_i+k_i\delta\in \Delta$. 
Then $\Pi_a\stackrel{\mbox{def}}{=}\bigcup_{i=1}^s \{\Pi_i,\beta_i\}$ is a set of simple roots of an affine root subsystem 
 of $\Delta$ that contains $\Delta_1$.

Reversing the procedure above, we obtain the following theorem.

\begin{theorem}
Let $\Delta$ be an affine root system and $\widetilde\Delta$ its affine root subsystem.
Let $\Delta_1,\dots,\Delta_k$ 
be indecomposable components of $\widetilde\Delta$. 
Let $\Pi$ be a set of simple roots of $\widetilde\Delta$ and $\Pi_i$
a set of simple roots of $\Delta_i$, $1\leq i\leq k$. 
Denote by $J\subset \{1,\dots,k\}$ the set of indices
such that $j\in J$ if $\Delta_j$ is a root system of an untwisted algebra.  Let $J'$ be a subset of $J$.
For each $j\in J'$ discard the root $\alpha_0^j\in \Pi_j$ from $\Pi$ (see \cite[Table Aff~1]{Kac})
and denote by $\Pi'$ the set of simple roots obtained. 
Let $\Delta'\subset \Delta$ be the root system whose set of simple roots is $\Pi'$. 

Then  $\Delta'$ is a root subsystem of  $\Delta$.

Moreover, for every non-affine root subsystem $\Delta'$ of $\Delta$, there exists an affine root subsystem $\widetilde\Delta$ of $\Delta$ such that $\Delta'$ can be obtained from $\widetilde\Delta$ in this way.
\end{theorem}

\section{Conformally invariant regular subalgebras of affine Kac--Moody algebras}\label{coset}

Recall that an untwisted affine Kac--Moody algebra is an extension of a loop algebra over a semisimple finite-dimensional Lie algebra $\fg$.  For the duration of this section, we fix the notation $\hat\fg=\fg[t,t^{-1}]\oplus\CC K$, where $K$ is a central element.  The bracket in $\hat\fg$ is defined as $[at^m,bt^n]=[a,b]t^{m+n}+\delta_{m,-n}(a|b)K$, where $(\ |\ )$ is an invariant form on $\fg$ such that the square of the highest root of $\fg$ is $2$.  Given an irreducible representation of $\hat\fg$, $K$ necessarily acts as a scalar $k$ called the level of this representation.  Though $k$ can take on any value, in physical literature it is assumed to be integral and positive.  For the most part we will impose some of these restrictions as well. However, it should be pointed out that from a purely mathematical standpoint, half-integer levels are also interesting  (see e.g. \cite{P} for a vertex algebra construction related to a pair of conformally invariant algebras discussed below).

It is well known that the vacuum module $V_k(\fg)$ of $\hat\fg$ of level $k$ carries the structure of a vertex operator algebra with the canonical conformal vector $\omega_k(\fg)$ given by the Sugawara constuction.  The central charge of $\omega_k(\fg)$ is $$D_k(\fg)=\sum \frac{k\dim\fg_i}{k+{\check{h}}_i},$$ where $\fg_i$ are irreducible components of $\fg$, $\fg=\oplus\fg_i$, and ${\check{h}}_i$ is the dual Coxeter number of $\fg_i$.  The Fourier coefficients of $\omega_k(\fg)$ form a Virasoro algebra that acts naturally on $\hat\fg$-modules of level $k$ \cite{Kac-va}. We can thus construct a large family of representations of the Virasoro algebra with central charge $D_k(\fg)$ with desired properties (unitary, of highest weight, etc).

\subsection{Coset construction} This construction can be generalized via the so-called coset construction as follows \cite{GKO1, Kac-va}.
Consider a graded subalgebra $\hat\fh\subset\hat\fg$ whose Cartan subalgebra is in the Cartan subalgebra of $\hat\fg$ and has the same imaginary root for some choice of root system.  (Apart from subalgebras of the same type, all regular subalgebras constructed above satisfy this condition.)  $V_{k'}(\fh)$ naturally embeds into $V_k(\fg)$ for a specific choice of $k'$.  Namely, for $\hat\fh$ simple, let $j$ be the square of the highest root in $\hat\fh$ with respect to the form $(\ |\ )$ on $\hat\fg$.  Then $k'=2k/j$.  (The quantity $2/j$ is known as the Dynkin index of $\fh$ in $\fg$).  In the semisimple case, $V_{k'}(\fh)$ is a direct sum of affine vertex algebras corresponding to simple components of $\fh$.
We can thus consider the difference $\omega_k(\fh,\fg)=\omega_k(\fg)-\omega_k(\fh)$.  A direct computation utilizing the Sugawara construction shows that $\omega_k(\fh,\fg)$ is conformal with the central charge $D_k(\fg,\fh)=D_k(\fg)-D_k(\fh)$.

If $D_k(\fg,\fh)\neq 0$, we obtain yet another construction of representations of the Virasoro algebra with a given charge.  In particular, this strategy was used in the seminal work \cite{GKO2} to construct all irreducible unitary highest weight representations of the Virasoro algebra with charges between 0 and 1.

The case $D_k(\fg)=D_k(\fh)$ is also important, since then the pair $\fh\subset\fg$ can be used to construct consistent conformally invariant theories.  Such pairs were studied in e.g. \cite{BB} or \cite{SW}.  Both papers imposed certain restrictions such as $k=1$ (or, more generally, positive integral) in \cite{BB} or algebras being simply-laced in \cite{SW}.  

Below, using our classification of regular subalgebras of simple non-twisted affine Kac--Moody algebras, we provide a unified strategy for computing all possible values of central charge for $\omega_k(\fh,\fg)$, where $\fh$ is a regular subalgebra of $\fg$.  We also indicate when $\fh$ is conformally invariant (i.e. when $D_k(\fg)=D_k(\fh)$).

Note first that if $\hat\fh$ is not a maximal subalgebra of $\hat\fg$, i.e. if $\hat\fh\subset\hat\fh'\subset\hat\fg$, then it suffices to compute the required values for the pairs $\hat\fh\subset\hat\fh'$ and $\hat\fh'\subset\hat\fg$.  Thus we restrict our attention to maximal subalgebras only.  As noted above, we will also consider only subalgebras of the type different from $\hat\fg$'s.  Moreover, if $\fg$ is not simple, i.e. $\fg=\oplus\fg_i$, then it suffices to consider projections of $\fh$ onto $\fg_i$ as subalgebras of $\fg_i$.  Thus, we can assume that $\fg$ is simple. 

\subsection{Conformally invariant subalgebras of non-twisted algebras} We consider each of the types of $\fg$ in turn.  The values of dual Coxeter numbers for $\hat\fg$ can be found in, e.g. \cite[6.1]{Kac}.  We assume that $k\neq \check{h}$ (i.e. that the level $k$ is not critical) for all algebras and subalgebras considered.

\medskip
\noindent 
$\mathbf A_n$.  A maximal regular subalgebra $\fh$ must be of the type $A_m+A_{n-m-1}$, $0<m<n$. (Here and below we freely use subalgebra lists in Theorems~\ref{th-ind}, \ref{th-dec}, and~\ref{lower}.)  Then $D_k(\fg,\fh)=k\dfrac{(n+2k+1)(n+k(n-m)(m+1))-m(n-m-1)+(k+1)^2}{(k+n+1)(k+m+1)(k+n-m)}.$  Clearly for $k,n,m>0$, $n+2k+1>m$ and $n+k(n-m)(m+1)>n-m-1$, hence the numerator in $D_k(\fg,\fh)$ is positive.  It follows that there are no conformally invariant maximal regular subalgebras for $k$ positive.  More generally, a regular subalgebra of $\fg$ of type $A_n$ is a direct sum of subalgebras of the same type.  Thus it suffices to consider a subalgebra $\fh$ of type $A_m$.  Here $D_k(\fg,\fh)=k\dfrac{(n-m)((n+m+2)(k+1)+nm)}{(k+n+1)(k+m+1)}$.  Thus there are no conformally invariant subalgebras in type $A_m$ for $k$ positive.

\medskip
\noindent 
$\mathbf B_n$.  For $\fh$ of type $B_{n-1}$, $D_k(\fg,\fh)=k\dfrac{4n^2-8n+4kn-k+1}{(k+2n-1)(k+2n-3)}$.  The numerator, if viewed as a function of $n$, has discriminant $k^2-3k+3$.  For an integer $k$, it is a full square only for $k=2$.  Thus the numerator is never zero if $k$ and $n$ are integral and there are no conformally invariant pairs $(\fg,\fh)$ of this type.

For $\fh$ of type $D_n$, $D_k(\fg,\fh)=kn\dfrac{2k+2n-3}{(k+2n-1)(k+2(n-1))}$.  Clearly, there are no conformally invariant pairs for any $n$ on every integer level.

For a decomposable subalgebra $\fh$ of type $B_m+D_{n-m}$, $0<m<n$, $D_k(\fg,\fh)=k(n-m)(2m+1)\dfrac{(k-1)(2n+2k-3)}{(k+2n-1)(k+2m-1)(k+2(n-m)-2)}$.  The subalgebra is conformally invariant on level $1$.  Note that just as in the case of subalgebra of type $D_n$, we also get a conformally invariant subalgebra on the same half-integer level $-n-3/2$.

\medskip
\noindent 
$\mathbf C_n$.  For $\fh$ of type $A_{n-1}$, $D_k(\fg,\fh)=k\dfrac{n^3+n^2k+2kn+k+n+1}{(k+n+1)(2k+n)}$.  Setting the numerator to zero and solving for integral $k$, we see that $n\equiv 1\mbox{\ mod\ }3$.  Also, a direct computation shows that $-n/3<k<1-n/3$.  Both statements imply that $k=n/3-1$ but for such value of $k$ the numerator is zero only when $n=-2$.  Hence there are no conformally invariant pairs of this type.

If $\fh$ is decomposable of type $C_m+C_{n-m}$, $0<m<n$, then $D_k(\fg,\fh)=m(n-m)k\dfrac{(2k+1)(2k+n+2)}{(k+n+1)(k+m+1)(k+n-m+1)}$.  Hence $\fh$ is conformally invariant only for even $n$ on the level $k=-n/2-1$.  Also, if $\fh$ is indecomposable of type $C_m$, $D_k(\fg,\fh)=k(n-m)\dfrac{(2(m+n)+1)(k+1)+2mn}{(k+n+1)(k+m+1)}$.  For $k$ positive, there are no conformally invariant subalgebras of this type.

\medskip
\noindent 
$\mathbf D_n$.  For $\fh$ of type $A_{n-1}$, $D_k(\fg,\fh)=k\dfrac{(k+1)n^2+(2-k)n+k-2}{(k+2(n-1))(k+n)}$.  The numerator, if viewed as a function of $n$, has discriminant $12-3k^2$.  For an integer $k$, it is a full square only for $k=1,2$, at which levels the numerator is not zero.  Thus $D_k(\fg,\fh)\neq 0$ for any integral $k$ and $n>1$ and there are no conformally invariant pairs of this type.

For a decomposable maximal subalgebra $\fh$ of type $D_m+D_{n-m}$, $0<m<n$, we have $D_k(\fg,\fh)=4km\dfrac{(k-1)(n-m)(n+k-2)}{(k+2(n-1))(k+2(m-1))(k+2(n-m-1))}.$  It follows that on level $k=1$ every subalgebra whose simple components are of type $D$ is conformally invariant.  Also, the maximal subalgebra is conformally invariant on level $2-k$.

\medskip
\noindent 
$\mathbf E_6$. $\fh$ is either indecomposable of type $D_5$ or decomposable of types $A_5+A_1$ or $A_2+A_2+A_2$.  Thus the possible values of $D_k(\fg,\fh)$ are $k\dfrac{84+33k}{(k+12)(k+8)}$, $40k\dfrac{(k-1)(k+3)}{(k+12)(k+6)(k+2)}$, $6k\dfrac{8k-21}{(k+12)(k+3)}$.  For $k=-3, 1$, the pair $(E_6, A_5+A_1)$ is conformally invariant.

\medskip
\noindent 
$\mathbf E_7$. $\fh$ can be of the types $E_6$, $A_7$, $D_6+A_1$, and $A_5+A_2$.  The corresponding values of $D_k(\fg,\fh)$ are $k\dfrac{55k+192}{(k+18)(k+12)}$, $70k\dfrac{k-1}{(k+8)(k+18)}$, $64k\dfrac{(k-1)(k+4)}{(k+18)(k+10)(k+2)}$, and $90k\dfrac{(k-1)(k+4)}{(k+18)(k+6)(k+3)}$.  Note that we obtain conformally invariant pairs for $k=1,-4$ and that subalgebras $D_6+A_1$ and $A_5+A_2$ produce conformally invariant pairs for the same levels.

\medskip
\noindent 
$\mathbf E_8$. $\fh$ is either indecomposable of the types $A_8$ and $D_8$ or decomposable of the types $E_7+A_1$, $E_6+A_2$, or $A_4+A_4$.  The possible values of $D_k(\fg,\fh)$ are $168k\dfrac{k-1}{(k+30)(k+9)}$, $128k\dfrac{k-1}{(k+30)(k+14)}$, $112k\dfrac{(k-1)(k+6)}{(k+30)(k+18)(k+2)}$, $162k\dfrac{(k-1)(k+6)}{(k+30)(k+12)(k+3)}$, and $200k\dfrac{k-1}{(k+30)(k+5)}$.  Note that all maximal subalgebras are conformally invariant for $k=1$ and decomposable subalgebras with nonisomorphic components are also both conformally invariant on level $k=-6$.  

\medskip
\noindent 
$\mathbf F_4$. $\fh$ can be of the types $B_4$, $A_2+A_2$, or $C_3+A_1$.  The possible values of $D_k(\fg,\fh)$ are $8k\dfrac{2k+5}{(k+9)(k+7)}$, $4k\dfrac{20k^2+51k+9}{(k+9)(k+3)(2k+3}$, and $k\dfrac{59k^2+61k-70}{2(k+9)(k+4)(k+1)}$.  There are no conformally invariant maximal regular subalgebras on integer levels (but there is a conformally invariant subalgebra on the half-integer level $-5/2$).

\medskip
\noindent 
$\mathbf G_2$. $\fh$ is either of the type $A_2$ or $A_1+A_1$.  Thus $D_k(\fg,\fh)$ is either $2k\dfrac{17k+5}{3(k+4)(k+1)}$ or $2k\dfrac{15k^2+26k+14}{(k+4)(k+2)(3k+2)}$.  There are no conformally invariant maximal regular subalgebras.
\medskip

\begin{remark} (1) Our computations confirm the conclusion of \cite{SW} that for positive integer levels conformally invariant subalgebras exist only on level $1$ (for a direct physical proof of this fact, see \cite{Ko}).  However, there are conformally invariant subalgebras on negative integer and half-integer levels.

(2) The above discussion implies that conformally
invariant subalgebras are rather rare.  However, every maximal 
subalgebra of an algebra of type $E_8$ is conformally invariant.
Such exceptional behavior is probably a consequence of some intrinsic
property of affine algebras of type $E_8^{(1)}$

\end{remark}

\subsection{Conformally invariant subalgebras of twisted affine algebras} There are analogs of the Sugawara and coset constructions of Virasoro algebras for twisted affine Kac--Moody algebras.  Details can be found in, e.g. \cite[3.4]{KW}.  For our purposes it suffices to desribe what central charges we can obtain from an algebra $\hat\fg(\sigma)$ of type $X_n^{(p)}$ and its subalgebra $\hat\fh(\sigma|_{\fh})$ of type $Y_m^{(q)}$.  Here $\sigma$ is an automorphism of $\fg$ that restricts to $\fh$.  On the level $k$, we obtain the realization of the Virasoro algebra with the charge $D_k(\fg)$, and corresponding to $\hat\fh(\sigma|_{\fh})$, another realization with the charge $D_{k'}(\fh)$.

Recall that $\sigma$ is completely determined by an $(m+1)$-tuple $(s_0,\dots,s_l)$, where $l$ is the rank of $\fg$: $\sigma(E_j)=\exp\left(2\pi is_j\left/\left(p\sum a_ks_k\right)\right.\right)E_j$.  Let $n_s(\fg)=\sum a_is_i$ (here $a_i$ are labels on the Dynkin diagram of $\hat\fg$).  For the pair $\hat\fh(\tau)\subset\hat\fg(\sigma)$ set 
$$r=\frac{qn_s(\fh)}{pn_s(\fg)}\quad\mbox{and}\quad D_k(\fg,\fh)=rD_k(\fg)-D_{k'}(\fh).$$
Below we will provide examples only for $\sigma$'s induced by diagram automorphisms.  Here $n_s(\fg)=n_s(\fh)=1$.

\medskip
\noindent
$\mathbf A_{2n}^{(2)}$. We consider the algebra $A_2^{(2)}$ separately.  It has two distinct maximal subalgebras of the type $A_1^{(1)}$.  The corresponding values of $D_k(\fg,\fh)$ are $k\dfrac{k-1}{(k+3)(k+2)}$ and $k\dfrac{13k-1}{2(k+3)(2k+1)}$.  Thus there is a conformally invariant subalgebra on level $1$.

For $n>1$, $A_{2n}^{(2)}$ has subalgerbas of the type $A_{2n-1}^{(2)}, B_n^{(1)}, A_{2m}^{(2)}+A_{2(n-m)-1}^{(2)}, D_m^{(1)}+A_{2n-2m}^{(2)}, D_3^{(2)}+A_{2n-4}^{(2)}$.

For a subalgebra of the type $A_{2n-1}^{(2)}$,  $D_k(\fg,fh)=k\dfrac{(4n^2+1)(k+1)+2n(4k+1)}{2(k+2n+1)(k+n)}$.  Direct computation shows that no integer $k$ makes this expression zero.  Thus a subalgebra of this type is never conformally invariant.
For $n=2$, the subalgebra obtained in the same manner has type $D_3^{(2)}$.  Here $D_k(\fg,fh)=5k\dfrac{29k+37}{2(k+9)(k+2)}$ and, again, the subalgebra is not conformally invariant.

For a subalgebra of the type $B_n^{(1)}$, $D_k(\fg,\fh)=kn\dfrac{k-2n-3}{(k+2n+1)(k+2n-1)}$.  The subalgebra is conformally invariant on the level $k=2n+3$.

For decomposable subalgebras of $A_{2n}^{(2)}$, we will write $D_k(\fg,\fh)$ in terms of $N=2n$ and $M=2m$. For the type $A_{2m}^{(2)}+A_{2(n-m)-1}^{(2)}$,
$D_k(\fg,\fh)=k[k(2(N-M)(M+N+MN)+
k(4(M+M^2-N)-(M+N)^2)+(5N-3M)) -(N-1)(M-1)+M^2+N^2+(k+1)^2+1]/\left[(k+N+1)(k+M+1)(2k+N-M)\right]$.  There are no conformally invariant subalgebras for small positive values of $k$.

For the type $D_m^{(1)}+A_{2n-2m}^{(2)}$, we have $$D_k(\fg,\fh)=km\dfrac{2Nk(N-M+2k+1)+M(N-k-3k^2+2)+5k^2+2k-2}{2(k+N+1)(k+N-M+1)(k+m-1)}.$$  Again, there are no conformally invariant subalgebras for small positive values of $k$.

For the type $D_3^{(2)}+A_{2n-4}^{(2)}$, we have $$D_k(\fg,\fh)=k\dfrac{2kN(8(k+N)-7)+(N-1)^2-31k^2-18k+12}{2(k+N+1)(k+N-3)(k+2)}.$$  The numerator is never zero mod 4, hence there are no conformally invariant subalgebras of this type.

\medskip
\noindent
$\mathbf A_{2n-1}^{(2)}$.  Maximal regular subalgebras have the type $A_{n-1}^{(1)}, C_n^{(1)}, A_{2m-1}^{(2)}+A_{2(n-m)-1}^{(2)}$.

For the type $A_{n-1}^{(1)}$, $D_k(\fg,\fh)=3nk\dfrac{2nk+1}{2(k+2n)(2k+n)}$.  There are no conformally invariant pairs of these types.

For the type $C_n^{(1)}$, $D_k(\fg,\fh)=-(2n+1)k\dfrac{4n^2-2n+2+(2n+1)k}{2(k+2n)(k+2n+2)}$.  This expression is zero only if $k\in(-2n-1,-2n-2)$, hence for integral $k$ there are no conformally invariant subalgebras of this type.

For the type $A_{2m-1}^{(2)}+A_{2(n-m)-1}^{(2)}$, $D_k(\fg,\fh)$ can be found from the expression for the subalgebra of type $A_m^{(1)}+A_{n-m-1}^{(1)}$ of the algebra $A_n^{(1)}$ above.  Here we also conclude that there are no conformally invariant subalgebras on positive levels.

\medskip
\noindent
$\mathbf D_{n+1}^{(2)}$.  For a generic $n$, maximal regular subalgebras have the type $B_n^{(1)}$, $D_n^{(2)}$, and $D_{n-m}^{(1)}+D_{m+1}^{(2)}$.  

For the type $B_n^{(1)}$, $D_k(\fg,\fh)=-k\dfrac{2(2n^2-n+1)+(3n-1)k}{2(k+2n)(k+4n-2)}$.  $D_k(\fg,\fh)=0$ for $k=-(4n^2-2n+2)/(3n-1)$, hence $k\in(-4n/3-7/9,-4n/3+2/9)$.  Depending on the value of $n\mbox{\ mod\ } 3$, $k$ must be of the form $4n/3, (4n+1)/3, (4n+2)/3)$.  Only in the former case and only for $n=3$, $k$ makes $D_k(\fg,\fh)$ zero.  Thus we have a conformally invariant subalgebra of $D_3^{(2)}$ on level $-4$.

For the type $D_n^{(2)}$, $D_k(\fg,\fh)=k\dfrac{4n(n+k-1)+k-2}{(k+2n)(k+2n-2)}$.  This expression is zero only if $k\in(n-2,n-1)$, hence there are no conformally invariant pairs of this type.

For the type $D_{n-m}^{(1)}+D_{m+1}^{(2)}$, $D_k(\fg,\fh)=-k(n-m)[8((n-m)(2mn+1+k(m+n))-1)+2k(3n(1-2k)+m(6-5k)-2mn)+7k(2-k)]/[2(k+2n)(k+2m)(k+4n-4m-4)]$.  We note that this subalgebra is not conformally invariant for $k=1$.  

An algebra of type $D_3^{(2)}$ also has a maximal subalgebra of type $C_2^{(1)}$.  For this pair, $D_k(\fg,\fh)=-5k\dfrac{5k+14}{2(k+4)(k+6)}$.

\medskip
\noindent
$\mathbf E_6^{(2)}$.  A maximal regular subalgebra can have the type $C_4^{(1)}, F_4^{(1)}, A_5^{(2)}+A_1^{(1)}$.  The respective values of $D_k(\fg,\fh)$ are $-3k\dfrac{11k+158}{(k+12)(k+10)}$, $-13k\dfrac{5k+42}{(k+12)(k+18)}$, and $k\dfrac{31k^2+4k-672}{2(k+12)(k+6)(k+4)}$.  Thus $C_4^{(1)}$ and $F_4^{(1)}$ subalgebra are never conformally invariant.  The subalgebra of the type $A_5^{(2)}+A_1^{(1)}$ is not conformally invariant as well; direct computations also show that neither are its simple components.

\medskip
\noindent
$\mathbf D_4^{(3)}$.  Maximal regualr subalgebras have the type $G_2^{(1)}, A_2^{(1)}, A_1^{(1)}+A_1^{(1)}$.  The corresponding  values of $D_k(\fg,\fh)$ are $-28k\dfrac{k+3}{(k+6)(k+12)}$, $6k\dfrac{k-1}{(k+6)(k+3)}$, and $2k\dfrac{k-4}{(k+6)(k+2)}$.  We obtain conformally invariant maximal subalgebras on the levels $k=-3, 1, 4$.

\section{Regular non-indefinite subalgebras of hyperbolic Kac--Moody algebras}
\label{hyp}

In this section we apply the classification of regular subalgebras of affine Kac--Moody algebras
to the investigation of regular subalgebras of hyperbolic Kac--Moody algebras.

Let $\Delta$ be a hyperbolic root system (i.e., a root system of a hyperbolic Kac--Moody algebra).  Denote by $W$ its Weyl group. Let $\Delta_1\subset\Delta$ be a root subsystem containing only finite and affine indecomposable components. Then $\Delta_1$ is a subsystem of some maximal affine or finite root subsystem $\Delta'$ of $\Delta$.  
Root systems of hyperbolic Kac--Moody algebras are described by those Dynkin diagrams whose proper subdiagrams are of either finite or affine type. Maximal subdiagrams correspond exactly to maximal affine or finite root subsystems.  Their Weyl groups can be obtained as groups generated by all but one reflections in the facets of the fundamental chamber of $W$. 

Given a hyperbolic Kac---Moody algebra $\fg$, the description above allows us to classify all regular subalgebras  
containing no summands of indefinite type. More precisely, any such subalgebra can be obtained by the following procedure. 

Take the Dynkin diagram $\Sigma$ of the root system of $\fg$ and denote the number of its nodes by $l$. Consider all  subdiagrams of $\Sigma$ with exactly $l-1$ nodes.  The corresponding subalgebras are either finite or affine.  Now take all regular subalgebras of the algebras obtained. For this, we either rely on the results of Section~\ref{af} (in case of affine algebras), or the results of~\cite{Dyn} (in case of semisimple finite-dimensional algebras).

\bigskip

\textbf{Acknowledgements.}\\
The first author was partially supported by grants NSh-5666.2006.1, INTAS YSF-06-10000014-5916, and  RFBR 07-01-00390-a.  The third author was partially supported by grants NSh-5666.2006.1, INTAS YSF-06-10000014-5766, and  RFBR 07-01-00390-a.

We thank the referees for valuable comments.

\end{document}